\documentclass[10pt]{amsart}
\usepackage{amsmath, latexsym}
\usepackage{amssymb,amsthm}
\usepackage{amscd, 
}
 
\newtheorem{theorem}{Theorem}[section]
\newtheorem{proposition}[theorem]{Proposition}
\newtheorem{lemma}[theorem]{Lemma}

\theoremstyle{remark}
\newtheorem{remark}[theorem]{Remark}
\theoremstyle{remark}
\newtheorem{exercise}[theorem]{Exercise}

\theoremstyle{definition}
\newtheorem*{definition}{Definition}

\def\C{\mathbb{C}}
\def\real{\mathbb{R}}
\def\integer{\mathbb{Z}}
\def\supp{\mathrm{supp}}
\def\id{\mathrm{id}}
\def\sphere{\mathbf{S}^{d-1}}
\def\cone{\mathbf{C}} 
\def\cB{\mathcal{B}}
\def\BB{\mathcal{C}}
\def\FF{\mathcal{F}}
\def\MM{\mathcal{ M}}

\def\TT{{\mathcal{T}}}

\newcommand{\hol}[1]{_{C_*^{#1}}}%
\def\LL{\mathcal{L}}
\def\L{\mathcal{L}}

\def\gg{{\gamma}}

\newcommand{\cout}[1]{} 

\newcommand{\comment}[1]{}

\begin{document}

\title{Spectra of differentiable hyperbolic maps}

\author{Viviane Baladi and Masato Tsujii}

\address{CNRS-UMR 7586, Institut de Math\'e\-ma\-ti\-ques  Jussieu, Paris, France}
\email{baladi@math.jussieu.fr}

\address{Mathematics, Hokkaido University, 
Sapporo, Hokkaido, Japan}

\email{tsujii@math.sci.hokudai.ac.jp}

\date{March 2006}

\begin{abstract}
This note is about  the spectral
properties of transfer operators  
associated to smooth hyperbolic
dynamics.
In the first two sections, we state our new results \cite{BT2}
relating such spectra with dynamical determinants, first
announced at the conference ``Traces in Geometry, Number Theory and Quantum Fields" 
at the Max Planck Institute, Bonn, October 2005.
In the last two sections,  we give a reader-friendly presentation of
some  key ideas in our work 
in the simplest possible settings, including a new proof
of a result of Ruelle on expanding endomorphisms.
(These last two sections are  a revised version of the
lecture notes  given during the
workshop ``Resonances and Periodic Orbits:
Spectrum and Zeta functions in Quantum and Classical Chaos"
at Institut Henri Poincar\'e, Paris, July 2005.)
\end{abstract}

\thanks{
This version has benefitted from the
remarks of S. Gou\"ezel and G. Keller. We thank them both warmly.
}

\maketitle


\section{A brief introduction}
For smooth  hyperbolic dynamical systems  and smooth weights
(smooth means $C^r$ for $r>1$),
we announce new results from \cite{BT2} relating  Ruelle transfer operators with
dynamical Fredholm determinants and dynamical zeta
functions:
First we establish bounds for the essential spectral radii 
of the transfer operator
on new spaces of anisotropic distributions
(Theorem~\ref{mainprop} and Lemma~\ref{Kitform}), 
improving previous results  (Theorem~\ref{mainprop0} from \cite{BT}), and
giving  variational expressions for the bounds.
Then (Theorem~\ref{main}),
we give a new proof of Kitaev's \cite{Ki} lower bound for the radius of  the disc
in which the dynamical Fredholm determinant admits
a holomorphic extension, and,
in addition, we show that  the zeroes of the determinant in the corresponding disc 
are in bijection with
the  eigenvalues of the   transfer operator on our spaces.
The proofs are based on elementary Paley-Littlewood analysis in Fourier
space, using (and improving) a decomposition of the Fourier space into stable and
unstable cones, inspired by \cite{AGT} and introduced in \cite{BT}. 
To prove the results on the dynamical determinants
we   introduce  in \cite{BT2} methods based
on approximation numbers \cite{Pie}.

In Section \ref{state}, we give precise definitions and statements
of our new results (proofs will appear
elsewhere \cite{BT2}), recalling also some previous
results from \cite{BT}. Sections \ref{toy} and \ref{S1}   contain a ---
hopefully --- pedagogical
presentation of several key ideas and techniques in the proofs in two
simple, but nontrivial, cases (a few steps of the argument are
left
as exercises for the reader):

In Section~\ref{toy} we discuss, as a warm-up,
transfer operators associated to smooth ($C^r$ for
$r>1$) expanding endomorphisms on a manifold $X$.
The case of expanding, noninvertible, maps is easier than the case
of hyperbolic, invertible, maps, because composition by each local inverse branch
improves regularity, and a relevant Banach space
is $C^{r}(X)$. The  bounds  on the essential spectral radius 
together with the connection with the dynamical determinants are well-known
(see \cite{Ruelle1}, \cite{Ru2}, \cite{GuL}).
We give a new proof of
of the bounds of the essential spectral radius (this proof is
the only  original material in this text). This allows us to recall the basic Paley-Littlewood (or dyadic)
decomposition tools that are instrumental in \cite{BT} and \cite{BT2}.

In Section~\ref{S1}, we consider the simplest hyperbolic diffeomorphisms,
Anosov maps, giving the definition
of a Banach space of distributions suitable
for the hyperbolic case,  and
explaining the key steps in the proof
of the bounds in \cite{BT} and \cite{BT2}
on the essential spectral radius of the transfer
operators.


\section{New results on transfer operators and dynamical determinants}
\label{state}

Let $X$  be a $d$-dimensional $C^\infty$ Riemann manifold,
and let $T:X \to X$ be  a  
diffeomorphism which is of class $C^r$ for some $r>1$. (If $r$ is not an integer, 
this means that the derivatives of $T$ of order $[r]$ satisfy an
$r-[r]$ H\"older condition.)  Assume that there exists
a hyperbolic basic set $\Lambda\subset X$ for $T$.
This means that $\Lambda$ is $T$-invariant, transitive 
and
that there exist a compact  
neighborhood $V$ of $\Lambda$ such that $\Lambda=\cap_{m\in \integer}T^{m}(V)$ and 
an invariant decomposition 
$T_{\Lambda}X=E^u\oplus E^s$ (with $E^u\ne 0$
and $E^s\ne 0$) of the tangent bundle over $\Lambda$, such that 
for some constants 
$C>0$ and $0<\lambda_s<1$, $\nu_u>1$,
we have for all $m\ge 0$ and $x\in \Lambda$
\begin{equation}\label{hypexp2}
\|DT^m|_{E^s}\|\le C\lambda_s^{m}\quad\mbox{and}
 \quad \|DT^{-m}|_{E^u}\|\le C\nu_u^{-m}\, .
\end{equation}

For $s\ge 0$, let $C^s(V)$ be the set of complex-valued $C^s$ functions on $X$ with support  
contained in the interior of $V$. 
The Ruelle transfer operator associated to the
dynamics $T$ and the weight $g\in C^{r-1}(V)$ is defined by
\[
\LL=\LL_{T,g}:C^{r-1}(V)\to C^{r-1}(V),\quad \LL \varphi(x)=
g(x)\cdot \varphi\circ T(x \, ).
\]
Since $T$ is  hyperbolic, $\LL$ is not smoothness improving, so that it is in fact
not very interesting to let $\LL$ act on spaces of smooth functions.
One of our goals is to find a space of distributions on $V$ which is not
too small (it should contain
all $C^{r-1}$ functions ) and not too large ($\LL$ should be bounded, with
some control on the 
\footnote{See \S~\ref{toy} for a definition of the essential spectral radius.}
essential spectral radius, guaranteeing in particular that $\LL$ is
quasicompact, and that $\LL$ has
a spectral gap when $g$ is strictly positive on $V$). In other words, we
are aiming at yet another avatar of the Ruelle-Perron-Frobenius theory
in infinite dimension.  (See e.g. \cite{Ba0} for more classical examples.)

Our latest result in this direction
improves the bounds of \cite{BT} and \cite{GL} on the
spectrum of $\LL$
(We refer to the introduction of \cite{BT}
for historical comments
and references to the previous works,
\cite{BKL}, \cite{Ba}, and in particular the important
paper of Gou\"ezel and Liverani \cite{GL}.)
To state it, we need some notation
(see \cite{W} for background on ergodic theory).
For a 
$T$-invariant Borel probability measure $\mu$ on $\Lambda$, we write
$h_\mu$ for the metric entropy of $(\mu,T)$, 
and $\chi_\mu(A)\in \real \cup \{-\infty\}$ for the largest Lyapunov
exponent of a linear cocycle $A$ over $T|_\Lambda$,
with $(\log \|A\|)^+\in L^1(d\mu)$.
Let $\MM(\Lambda,T)$ denote the set of  $T$-invariant 
ergodic Borel probability
measures  on  $\Lambda$.

\begin{theorem}[Bounds on the essential spectral radius
\cite{BT2}]\label{mainprop}
Let $r>1$, $T$,  and $\Lambda \subset V$ be as above.
For any real numbers $q<0<p$  so that $p-q<r-1$,
there exists a Banach space $\BB^{p,q}(T,V)$ of distributions on $V$,
containing $C^s(V)$ for any $s>p$, and  contained
in the dual space of $C^{s}(V)$ for any $s>|q|$, with the following
properties:

For any $g\in C^{r-1}(V)$, the Ruelle operator $\LL_{T,g}$ 
extends to a bounded operator on $\BB^{p,q}(T,V)$. Its
essential spectral radius on this space is  not larger than
\begin{align*}
&Q^{p,q}(T,g)=\\
\nonumber&\,
\exp \sup_{\mu \in \MM(\Lambda,T)}
\Bigl \{h_\mu + \chi_\mu\left (\frac{g}{\det (DT|_{E^u})}\right ) +
\max\bigl \{p \chi_\mu(DT|_{E^s}), |q| \chi_\mu(DT^{-1}|_{E^u} )\bigr \}
\Bigr \}\, .
\end{align*}
\end{theorem}

See \cite[\S 8]{CL}--\cite{GuL} for a variational expression
analogous to  $Q^{p,q}(T,g)$
in the setting of $C^r$ expanding endomorphisms. 
Note that $\chi_\mu(g/\det (DT|_{E^u}))=\int \log |g|\, d\mu
-\int \log |\det( DT|_{E^u})|\, d\mu$, but the 
expression  as a Lyapunov exponent
is useful when $g$ is replaced by a bundle automorphism
(see \cite{BT2}).

\begin{remark}[Decay of
correlations] Assume for a moment that
$\Lambda$ is attracting for $T$, i.e., $T(V)\subset \mbox{interior} (V)$.
Once we have the estimates in Theorem \ref{mainprop}, 
it is not difficult to see that the spectral radius of the pull-back operator 
$T^* \varphi=\varphi \circ T$ on  $\BB^{p,q}_*(T,V)$ is equal to one.
(The constant function is a fixed function.)
If $(T,\Lambda)$ is in addition  topologically
mixing, then $1$ is the unique eigenvalue on the unit circle, 
it is a simple eigenvalue, and
the fixed vector  of the dual
operator to $T^*$ gives rise to  the SRB measure $\mu$: This 
corresponds to exponential decay of correlations for $C^p$ observables and $\mu$. (See Blank--Keller--Liverani \cite[\S3.2]{BKL} for example.)
\end{remark}

\begin{remark} [Spectral stability] It is not difficult to see that
there is $\epsilon >0$  so that
if $\widetilde T$ and $\widetilde g$, respectively, are
$\epsilon$-close to $T$ and $g$, respectively, in the $C^r$, resp.
$C^{r-1}$, topology,  then the associated
operator $\L_{\widetilde T,\widetilde g}$
has the same spectral properties than $\L_{T,g}$
on {\it the same Banach spaces.} Spectral
stability can then be proved, as it has been done in \cite{BKL} or \cite{GL}
for the norms defined there. We refer to \cite{BT2} for details.
\end{remark}

\smallskip
We 
next give an alternative expression for $Q^{p,q}(T, g)$. 
If $g\in C^{0}(V)$, we write 
$$ 
g^{(m)}(x)=\prod_{k=0}^{m-1} g(T^k(x)) \,, \qquad \forall m\in \integer_+\, .
$$
Put $\lambda=\max\{\lambda_s,\nu_u^{-1}\}$.
We define local hyperbolicity exponents for $x\in \Lambda$ and $m\in \integer_+$
by
\begin{equation}
\begin{aligned}
\lambda_{x}(T^{m})&=\sup_{v\in E^s(x)\setminus\{0\}} 
\frac{\|DT^{m}_x(v)\|}{\|v\|} \le C\lambda^m \quad \mbox{and}\\
 \nu_x(T^{m})
&=\inf_{v\in E^u(x)\setminus\{0\}} 
\frac{\|DT^m_x(v)\|}{\|v\|}\ge C^{-1}\lambda^{-m}.
\end{aligned}
\end{equation}
For  arbitrary real numbers $q$, $p$  and integer $m\ge 1$, set for $x\in \Lambda$
\begin{align}\label{pqm}
&\lambda^{(p,q,m)}(x)=
\max\bigl \{ (\lambda_{x} (T^m))^{p}, (\nu_{x} (T^m))^{q}
\bigr \} \, .
\end{align}

We may extend $E^s(x)$ and $E^u(x)$ to continuous bundles on $V$
(which are not  invariant in general), so that the inequalities (\ref{hypexp2}) hold for $x\in \cap_{k=0}^{m-1}T^{-k}(V)$, and for all $m\ge 0$, with some constant $C$. We
may thus extend the definition of $\lambda_{x}(T^{m})$, $\nu_x(T^{m})$ and $\lambda^{(p,q,m)}(x)$ to $\cap_{k=0}^{m-1}T^{-k}(V)$.
Letting $dx$ denote Lebesgue measure on
$X$, define  for integers $m\ge 1$, and  $p$, $q\in \real$
\begin{equation}
\rho^{p,q}(T, g,m)= 
\int_{X}  |g^{(m)}(x)| \lambda^{(p,q,m)}(x)\,  dx \, .
\end{equation}
In  \cite{BT2}, we also show:

\begin{lemma} \label{Kitform}
For $r>1$, $T$,  and $V$ as above, and  $g\in C^\delta(V)$
for some $\delta >0$, the limit 
$\rho^{p,q}(T, g)=\lim_{m\to \infty}(\rho^{p,q}(T, g, m))^{1/m}$ 
exists for all $p$, $q\in \real$.
If $q\le 0\le p$, we have
$Q^{p,q}(T,g)= \rho^{p,q}(T, g)$.
\end{lemma}

Kitaev \cite{Ki} proved existence of the limit $\rho^{p,q}(T, g)$, and showed
that it gave a lower bound for the domain of holomorphic extension
of a dynamical determinant (see also our Theorem~\ref{main} below).

Next, we compare Theorem~\ref{mainprop}
to our previous results, using the above lemma.
(It is convenient to put $a^{1/\infty}=1$ for $a\in \real_+^*$.)
In \cite{BT}, we proved:

\begin{theorem}\label{mainprop0}
Let $r>1$, $T$,  and $\Lambda \subset V$ be as above.
For any real numbers $q<0<p$  so that $p-q<r-1$,
there exist  a Banach space  $\BB_*^{p,q}(T,V)$
of distributions on $V$, and for
each $1<t<\infty$,  a Banach space $W^{p,q,t}(T,V)$ of distributions on $V$, 
with the following
properties:

$\BB_*^{p,q}(T,V)$ and $W^{p,q,t}(T,V)$ 
both contain $C^s(V)$ for any $s>p$, and  they both are  contained
in the dual space of $C^{s}(V)$ for any $s>|q|$.

For any  $g\in C^{r-1}(V)$, the operator $\LL_{T,g}$ 
extends boundedly to  $\BB_*^{p,q}(T,V)$, with
essential spectral radius not larger than
$R^{p,q,\infty}(T,g)$, and $\LL_{T,g}$ extends
boundedly to
$W^{p,q,t}(V)$
with
essential spectral radius not larger than
$R^{p,q,t}(T,g)$, with
\begin{align*}
R^{p,q,t}(T,g)=\lim_{m\to \infty}\left(
 \sup_{\Lambda} |\det DT^m|^{1/t} |g^{(m)}(x)| \lambda^{(p,q,m)}(x)\right)^{1/m}\, .
\end{align*}
\end{theorem}

Since $\rho^{p,q}(T,g)\le \inf_{1<t\le \infty} R^{p, q,t}(T,g)$, and
the inequality
can be strict,  Theorem~\ref{mainprop} can be viewed as
an improvement of Theorem~\ref{mainprop0}.
Note however that the anisotropic
Sobolev spaces $W^{p,q,t}(T,V)$  have applications 
to situations with less hyperbolicity, 
such as  skew products  \cite{AGT} or  time-one
maps of expanding semi-flows \cite{T}.

\bigskip

We next turn to  dynamical Fredholm determinants. 
The dynamical Fredholm determinant $d_{\LL}(z)$ corresponding to the Ruelle 
transfer operator $\LL=\LL_{T,g}$ is
\begin{equation}\label{eqn:zeta}
d_{\LL}(z)=
\exp \left(- \sum_{m=1}^\infty\frac {z^m}{m}
\sum_{T^m(x)=x}
\frac{g^{(m)}(x)} {|\det (1-DT^m(x))|} \right)\,  .
\end{equation}
The power series in $z$ which is exponentiated
converges only if $|z|$ is sufficiently small. 
The main new result in \cite{BT2} is about the analytic continuation of 
$d_{\LL}(z)$:
\begin{theorem}\label{main}
Let $r>1$, $T$, $V$, and $g\in C^{r-1}(V)$ be as above.

\noindent The function $d_{\LL}(z)$ extends  holomorphically
to the disc of radius $(\rho_r(T,g))^{-1}$ with
$$\rho_r(T,g)=\inf_{q<0<p, \, p-q<r-1} \rho^{p,q}(T,g)\, .
$$

\noindent For any real numbers $q<0<p$  so that $p-q<r-1$,
and each $z$ with $|z|< (\rho^{p,q}(T, g))^{-1}$, 
we have $d_\LL(z)=0$  
if and only if $1/z$ is an  eigenvalue of $\LL$ on $\BB^{p,q}(T,V)$, 
and the order
of the zero coincides with the algebraic multiplicity of the eigenvalue.
\end{theorem}

\begin{remark}
The proof in \cite{BT2} implies  that
for any real numbers $q<0<p$  so that $p-q<r-1$,  each
$1<t<\infty$,
and each $z$ with $|z|< (R^{p,q,t}(T, g))^{-1}$, 
respectively $|z|< (R^{p,q,\infty}(T, g))^{-1}$,
we have $d_\LL(z)=0$  
if and only if $1/z$ is an  eigenvalue of $\LL$ on $W^{p,q,t}(T,V)$, 
respectively $\BB_*^{p,q}(T,V)$,
and the order
of the zero coincides with the algebraic multiplicity of the eigenvalue.
\end{remark}

Note that for analytic hyperbolic diffeomorphisms and weights, it has been
known for 30 years  that $d_\LL(z)$ is an entire function
when the dynamical foliations are analytic \cite{Ru}. More
recently, Rugh and Fried
\cite{Rugh, Fr} studied $d_\LL(z)$ in this analytic framework, without any assumption on the
foliations, giving a spectral interpretation of its zeroes.

In the case of finite differentiability $r$, the connection
between transfer operators
and dynamical determinants  has been well  understood
in the easier setting of expanding
endomorphisms since 15 years ago (see \cite{Ru2}).
The case of hyperbolic diffeomorphisms 
has only been studied recently. 
In an important and pioneering article \cite{Ki}, Kitaev 
obtained the first claim of our Theorem ~\ref{main}, without the spectral interpretation
of the zeroes of $d_\LL(z)$.
Our proof is different and gives the spectral interpretation
of the zeroes of $d_\LL(z)$ contained in the second claim of  Theorem ~\ref{main}.
Note that a spectral interpretation of the zeroes
(in a smaller disc, depending on the dimension $d$)
has been obtained previously by Liverani \cite{Li}, 
using Banach spaces of \cite{GL}.

We refer to \cite{BT2} for the proof of Theorem~\ref{main}.


\section{A toy model: expanding endomorphisms}
\label{toy}

In order to give in the next section
the key ideas in the proof of Theorem~\ref{mainprop0}
in \cite{BT}
and 
Theorems~\ref{mainprop} and  ~\ref{main} in \cite{BT2}, we  revisit in this section
the much easier (and well understood) situation of expanding endomorphisms.
We first recall a definition and a few elementary facts, which will also be used
in Section~\ref{S1}. 

\begin{definition} [Essential spectral radius]
The essential spectral radius $r_{ess}(\L|_\cB)$ of a bounded operator
$\L$ on a Banach space $\cB$ is the infimum of the real numbers $\rho>0$
so that, outside of the disc of radius $\rho$, the spectrum of $\L$ on $\cB$
consists of isolated eigenvalues of finite multiplicity.
\end{definition}

The following basic fact will be at the very center of our proof
(it is behind  most techniques to estimate the essential
spectral radius: Lasota-Yorke or Doeblin-Fortet bounds, Hennion's theorem,
the Nussbaum formula, see e.g. \cite{Ba0}):

{\bf Compact perturbation.}
If  $\L=\L_1+\L_0$
where $\L_1$ is compact on $\cB$ and $\L_0$ is bounded on $\cB$, then
the essential spectral radius of $\L$ acting on $\cB$
is  not larger than the spectral radius of $\L_0$ acting on $\cB$.  
(See e.g. \cite{DS}.)

\smallskip

Not surprisingly, our main tool is integration by parts:

{\bf Integration by parts.}
By ``integration by parts on $w$," we will mean application, 
for $f\in C^2(\real^d)$ with
$\sum_{j=1}^{d}(\partial_j f(w))^2\ne 0$
and a compactly supported $g\in C^1(\real^d)$,
of the formula
\begin{align}\label{parts}
 \int e^{if(w)}g(w) dw
&=
-\sum_{k=1}^d 
\int i(\partial_k f(w))e^{if(w)}\cdot 
\frac{i(\partial_k f(w))\cdot g(w)}
{\sum_{j=1}^{d}(\partial_j f(w))^2} dw
\\
&=
i\cdot \int e^{if(w)}\cdot \sum_{k=1}^d \partial_k\left(\frac
{\partial_k f(w)\cdot  g(w)}
{\sum_{j=1}^{d}(\partial_j f(w))^2}\right) dw \, ,\notag
\end{align}
where $w=(w_k)_{k=1}^{d}\in \real^d$, and $\partial_k$ denotes partial differentiation 
with respect to $w_k$.  (Note that if $f$ is $C^ r$ we
can only integrate by parts $[r]-1$ times in the above sense, even
if $g$ is $C^ r$ and compactly supported.)

\smallskip
{\bf Regularised integration by parts}
If  $f\in C^{1+\delta}(\real^d)$ and  $g\in C^\delta_0(\real^d)$,
for $\delta\in (0,1)$, and $\sum_{j=1}^{d}(\partial_j f)^2\ne 0$ on $\supp(g)$,
we shall consider the following
``regularised integration by parts:" 
Set, for
$k=1,\ldots, d$
\begin{align*}
h_k:=
\frac{i(\partial_k f(w))\cdot g(w)}
{\sum_{j=1}^{d}(\partial_j f(w))^2} .
\end{align*}
Each $h_k$ belongs to $C^\delta_0(\real^d)$.
Let $h_{k,\epsilon}$, for small $\epsilon >0$, be the convolution of
$h_k$ with $\epsilon^{-d} \upsilon(x/\epsilon)$, where 
the $C^\infty$ function
$\upsilon :\real^d \to \real_+$ is supported in the unit
ball and satisfies $\int \upsilon(x)dx=1$. 
There is $C$, independent of $f$ and $g$,
so that for each small $\epsilon>0$ and all $k$,
$$
\| \partial_k h_{k,\epsilon} \|_{L^\infty}
\le C \|h_k\|_{C^\delta} \epsilon^{\delta-1} ,
\quad
\|  h_k - h_{k,\epsilon} \|_{L^\infty}
\le C \|h_k\|_{C^\delta}  \epsilon^{\delta} .
$$
Finally, for every real number $\Lambda \ge 1$
\begin{align}\label{regparts}
 \int e^{i\Lambda f(w)}g(w) dw
&=-\sum_{k=1}^d \int i\partial_k f(w) e^{i\Lambda f(w)}\cdot  h_k(w) dw\\
\nonumber &= \int \frac{e^{i\Lambda f(w)}}
{\Lambda} \cdot \sum_{k=1}^d \partial_k h_{k,\epsilon} (w) dw \\
\nonumber &\quad-\sum_{k=1}^d 
\int i\partial_k f(w) e^{i\Lambda f(w)}\cdot 
(h_k(w)-h_{k,\epsilon} (w)) dw
.
\end{align}


\subsection{The result for locally expanding maps}
\label{S0}

Let $T:X\to X$ be $C^r$ for $r>1$, where  $X$ is a $d$-dimensional
compact manifold.  In this section, we assume that $T$ is
a locally expanding map, i.e., there are $C>0$ and $\lambda_s <1$ so that for each
$x$, all $m\ge 1$ and all 
$v \in T_x X$, we have $\|D_x T^m v\| \ge C\lambda_s^ {-m} \|v\|$. The function $g$ is
assumed to be $C^r$. We study the operator 
$$
\L_{T^{-1},g} u (x):=\sum_{y: T(y)=x} g(y) u(y) \, .
$$
(This is
the transfer operator associated to the branches of $T^ {-1}$, which
contract by at least $\lambda_s$.)
Note that 
$$R(T^{-1},g):=\lim_{m\to \infty}
\Bigl ( \sup_x \sum_{y: T^ m(y)=x}
|g^ {(m)}(x)| \Bigr )^ {1/m}$$
 is the spectral radius
of $\L_{T^{-1},g}$ acting on continuous functions.
\footnote{If $g=|\det D T|^{-1}$ then it is well-known that
$R(T^{-1},g)=1$.}

For $p>0$, recall that the $C^p$ norm of
$u \in C^ \infty(\real^d)$ is 
\[
\|u\|_{C^p}=\max\left\{
\max_{|\alpha|\le [p]}
\sup_{x\in \real^d}
|\partial^\alpha u(x)|,
\max_{|\alpha|= [p]}
\sup_{x\in \real^d}\sup_{y\in \real^d/\{0\}} 
\frac{|\partial^\alpha u(x+y)-\partial^\alpha u(x)|}{\|y\|^{p-[p]}}
\right\}
\]
where $\partial^\alpha u$ for a multi-index
$\alpha=(\alpha_1, \ldots ,\alpha_d)\in \integer_+$ denotes
the partial derivative $\partial^{\alpha_1}_1 \cdots \partial^{\alpha_d}_d u$,
and $|\alpha|=\sum_{j=1}^d \alpha_j$.
For $\varphi \in C^p(X)$, the
above norm can be used in charts to define a norm $\|\varphi\|_{C^p(X)}$.

We shall prove  the following  result:

\begin{theorem}[Essential spectral radius for expanding maps]\label{th:exp}
Let $T$ be $C^r$ and expanding, and let $g$ be $C^{r}$ for $r>1$.
For any noninteger $0< p \le r$,  let $C^p_*(X)$ be the closure of $C^ 
\infty(X)$
for the $C^p$ norm. Then the operator $\L_{T^{-1}, g}$ is bounded on $C^p_*(X)$
and
$$r_{ess}(\L_{T^{-1},g}|_{C^ p_*(X)}) \le R(T^{-1},g)\cdot \lambda_s^ 
{p}\, 
.
$$
\end{theorem}

The main interest of the proof given here 
is that it can be generalised to the hyperbolic
case. Note also that if $p$ is an integer, the proof below gives the same bounds for a Zygmund space
$C^p_*(X)$.

\begin{exercise}
Prove that for any noninteger $p>0$ we have
$C^p_*(X) \subset C^p(X)$, and that the inclusion is strict.
\end{exercise}

Ruelle \cite{Ruelle1} proved the statement of Theorem~\ref{th:exp}
for $C^p(X)$ instead of $C^p_*(X)$.
(It is in fact possible to modify the definitions in Subsection~\ref{deff}, 
to get a new proof of Ruelle's result. This modification is  cumbersome
when dealing with distributions in the later sections, and we do not
present it here. Adapting the argument in \cite{BT2} to the case of expanding
endomorphisms in the spirit of this section, it is even possible to recover the
optimal bounds in \cite{GuL}.)


\subsection{Local definition of H\"older norms in Fourier coordinates}
\label{deff}

We present here the ``dyadic decomposition" approach to
compactly supported
H\"older functions in $\real^d$ (for $d\ge 1$).
Fix a $C^\infty$ function $\chi:\real_+\to [0,1]$ with
\[
\chi(s)=
1, \quad \mbox{for $s\le 1$,}\qquad
\chi(s)=0, \quad \mbox{for $s\ge 2$.}
\] 
Define $\psi_n:\real^d\to [0,1]$ for $n\in \mathbb Z_+$, by
$\psi_0(\xi)=\chi(\|\xi\|)$, and
\begin{align*}
\psi_n(\xi)=\chi(2^{-n}\|\xi\|)-\chi(2^{-n+1}\|\xi\|) \, , \quad n\ge 1 \, 
.
\end{align*}
We have  
$1=\sum_{n= 0}^{\infty}\psi_{n}(\xi)$, and   $\supp(\psi_n)\subset 
\{\xi\mid 2^{n-1}\le \|\xi\|\le  2^{n+1}\}$ for $n\ge 1$. 
Also $\psi_n(\xi)=\psi_1(2^{-n+1}\xi)$ for $n\ge 1$. Thus,
for every multi-index $\alpha$, there exists a constant $C_\alpha$ such that  
$\|\partial^\alpha \psi_n\|_{L^\infty} \le C_\alpha 2^{-n|\alpha|}$ for all $n\ge 0$, 
and the inverse 
Fourier transform  of $\psi_n$,
$$\widehat \psi_n(x)=
(2\pi)^{-d} \int_{\real^d} 
 e^{ix\xi} \psi_n(\xi)  d\xi\, ,\quad x \in 
\real^d \, ,
$$  
decays rapidly in the sense of Schwartz. Furthermore we have 
$$\widehat \psi_n(x)=2^{d(n-1)} \widehat \psi_1(2^{n-1}x)$$
for $n\ge1$ and all $x$, and 
\begin{equation}\label{bded1}
\sup_n \int_{\real^d} |\widehat \psi_n(x)| dx < \infty \, .
\end{equation}

\begin{exercise} Prove the above claims on $\psi_n$ and $\widehat \psi_n$.
\end{exercise}

Fix a compact subset $K\subset \real^d$ with non-empty interior
and let  $C^{\infty}(K)$ be the space of complex-valued $C^\infty$ 
functions on $\real^d$ supported 
on $K$.  
Decompose each $u \in C^\infty(K)$  as
$u=\sum_{n\ge 0}u_n$, by defining for 
$n\in \mathbb Z_+$ and $x\in \real^d$
\begin{equation}\label{decomp}
u_n(x)=\psi_n(D) u(x):=(2\pi)^{-d}\int_K \int_{\real^d}  e^{i(x-y)\xi}\psi_n(\xi)u(y) dy d\xi \,  .
\end{equation}
Note that $u_n$ is not necessarily  supported in $K$, although
it satisfies good decay properties when $\|x\|\to \infty$: we say
that the operator
$\psi_n(D)$ is not a ``local" operator, but it is ``pseudo-local."
(See \cite{BT}. The pseudo-local estimates there are useful e.g.
to show the compactness results in  Propositions~\ref{cpct0}
and \ref{prop:cpt}.)

\begin{remark}
The notation $a(D)$ for the operator sending a compactly
supported $u\in C^\infty(\real^d)$ to
$$a(D) u(x):=(2\pi)^{-d}\int_K \int_{\real^d}  e^{i(x-y)\xi}a(\xi)u(y) dy d\xi
=(\widehat a * u)(x) \, ,$$
associated to  $a\in C^\infty(\real^d)$
so that $\partial^\alpha a (\xi)\le C_\alpha(a) \|\xi\|^{-|\alpha|}$
for each multi-index $\alpha$, stands for the ``pseudo-differential operator
associated to the symbol" $a$.  We shall not need any knowledge about
pseudodifferential operators, and shall not require symbols depending
on both $x$ and $\xi$.
\end{remark}

\begin{definition}[Little h\"older space $C^p_*(K)$]
For a real number $p>0$,   define on $C^{\infty}(K)$ the norm
\[
\|u\|\hol{p}=  \sup_{n\ge 0} \; 2^{pn} \|u_n \|_{L^\infty(\real^d)}  \, .
\]
The space  $C^{p}_*(K)$ is
the completion of $C^{\infty}(K)$ with respect 
to   $\|\cdot \|\hol{p}$.
\end{definition}

\begin{remark} It is known that if $p$
is not an integer then the norm $\|u\|\hol{p}$ is equivalent to the $C^p$ norm.
(See \cite[Appendix A]{Ta}.)
\end{remark} 

We shall not give a proof of the following, very standard, result (the proof
is based on the Ascoli-Arzel\`a lemma; see Proposition \ref{prop:cpt}
for an anisotropic analogue):

\begin{proposition}[Compact embeddings]\label{cpct0}
If $0<p'<p$ the inclusion $C_*^p(K) \subset C_*^ {p'}(K)$ is compact.
\end{proposition}

\subsection{Compact approximation for local maps}
\label{comcon}

Let $r >1$.  Let $K, K'\subset \real^d$ be compact subsets with 
non-empty interiors, 
and take 
a compact  neighbourhood $W$ 
of $K$.  Let $\TT:W\to K'$ be a $C^{r}$ diffeomorphism onto
its image (the reader should think of $\TT$ as being a local inverse branch
of an expanding map $T$, in charts).
Let $\gg:\real^d\to \C$ be a $C^{r-1}$ function 
supported in the interior of $K$.  
In this section we study a local transfer operator:
\[
L: C^{r-1}(K')\to C^{r-1}(K), 
\qquad L u(x)= \gg(x)\cdot u\circ \TT(x) \, .
\]

We define a ``weakest contraction\footnote{This is because we will consider contracting maps $\TT$ in the application to Theorem~\ref{th:exp}.}" exponent
\[
\|\TT\|_+=\sup_{x\in K}
\sup_{\xi\ne 0} \frac{\|D\TT_x^{tr}(\xi)\|}
{\|\xi\|}\, .
\]

The following result is the key to the proof of Theorem~\ref{th:exp}:

\begin{theorem}\label{th:main1} 
For any real number $p>0$ such that $p\le  r-1$,
and every compact $K_0$
contained in the interior
of $K$, there is  a constant $C$, so that for each  $C^ r$ map $\TT$ 
as above
and every $\gg$ in $C^ {r-1}(K_0)$
there is a compact operator
$L_1:C^ p_*(K')\to C^p_*(K)$  such that for any $u\in C_*^{p}(K')$
\begin{equation*}
\|L u-L_1 u \|_{C_*^{p}}\le C\|\gg\|_{L^\infty}\cdot 
 \|\TT\|_+^{p} \|u\|_{C^p_*} \, .
\end{equation*}
If $\gg\in C^r(K_0)$ then the condition on $p$ may be relaxed to $0< p \le  r$. 
\end{theorem}

We sketch how to deduce Theorem~\ref{th:exp} from Theorem~\ref{th:main1}: 
Take a system of local charts $\kappa_i:V_i\to K_i\subset
\real^d$, $1\le i\le k$, and a $C^\infty$ partition of unity $\phi_i:X\to [0,1]$
subordinate to the covering by $V_i$, that is, the support
of $\phi_i$ is contained in the interior of $V_i$.
(Then, $C^p_*(X)$ is embedded in the direct sum of the local $C^p_*(K_i)$ spaces.)  Consider an iterate $\LL^m$ of $\LL$ and define the operators $\LL_{ij}:C^{r-1}_*(X)\to C^{r-1}_*(X)$ by 
$\LL_{ij}^m\varphi(x)=\phi_j\cdot \LL^m (\phi_i \varphi)$ so that
$\LL^m=\sum_{i,j}\LL^m_{ij}$. Since the operator $\LL^m_{ij}$ may be viewed as an operator in local charts, we may apply Theorem~\ref{th:main1} to $\LL^m_{ij}$,
taking  $\TT$  to be a branch of the inverse of $\kappa_i\circ T^m\circ \kappa_j^{-1}$, and taking $\gg$ to be $(\phi_j\cdot \phi_i\circ T^m \cdot g^{(m)})\circ \kappa_j^{-1}\circ \TT$. Then we get that the essential spectral radius of the operator $\LL^m=\sum_{i,j}\LL^m_{ij}$ is bounded by 
\[
R_m:=C \cdot \lambda_s^{mp} \cdot \biggl (\sup_{x\in X}\sum_{y:T^m(y)=x}|g^{(m)}(y)|\biggr )\, ,
\]
for some constant $C$ independent of $m$ and $g$. (It is crucial that the constant
$C$ in Theorem~\ref{th:main1} is independent of $\TT$ and thus of
the iterate
$m$.)
Thus the essential spectral radius of $\LL$ is bounded by $(R_m)^{1/m}$. Considering large $m$, we obtain Theorem~\ref{th:exp}.

\begin{proof}[Proof of Theorem \ref{th:main1}] 
We need a couple more notations.
Recall the function $\chi$ from Section~\ref{deff}. Define
$\tilde \psi_{\ell}:\real^d \to [0,1]$ by 
\[
\tilde \psi_{\ell}(\xi)=
\begin{cases}\chi(2^{-\ell-1}\|\xi\|)-\chi(2^{-\ell+2}\|\xi\|),& \mbox{ if 
$\ell\ge 1$,}\\
\chi(2^{-1}\|\xi\|),& \mbox{ if $\ell=0$.}
\end{cases}
\]
Note that $\tilde \psi_{\ell}(\xi)=1$ 
if $\xi \in \supp (\psi_{\ell})$.

We write\footnote{By definition, if
$\ell \not\hookrightarrow n$ then $n > \ell - n(\TT)$ for some
$n(\TT)$ depending only on $\TT$. This feature will not be present
in the hyperbolic case.} 
\begin{align*}
&\bullet \ell \hookrightarrow n \hbox{ if  } 2^ n\le  \| \TT\|_+ 2^{\ell+4}\, , \\
&\bullet \ell \not\hookrightarrow n \hbox{ otherwise.}
\end{align*} 
By the definition of $\not\hookrightarrow$
there exists  an integer $N(\TT)>0$ such that 
\begin{equation}\label{lowerbd0}
\inf_x d(\supp(\psi_{n}), D\TT_{x}^{tr}(\supp(\tilde \psi_{\ell})))
\ge 2^{\max\{n,\ell\}-N(\TT)} 
\quad \mbox{if $\ell \not\hookrightarrow n$.}
\end{equation}

Let  $h:\real^d \to [0,1]$ be  a $C^\infty$ function
 supported in $K$ and $\equiv 1$ on $\supp(\gg)$.
Noting that $L (f)$ is well-defined if $f \in C^\infty(\real^d)$ 
because $\gg$ is supported in $K$,  we may define  $L_1$ and $L_0$ by
$L_j (f)=(M \circ L_j')(f)$ with $Mf=h\cdot f$, and
$L'_j u=\sum_n (L'_j u)_{(n)}$ with 
\[
(L'_0 u)_{(n)}
=\sum_{\ell: \ell\hookrightarrow n}
\psi_{n}(D) (L \, u_\ell)\, ,
\]
and
\[
(L'_1 u)_{(n)}=
\sum_{\ell: \ell \not\hookrightarrow n}
\psi_{n}(D) (L\,\tilde{\psi}_{\ell}(D) u_\ell) \, .
\]
Since $\tilde{\psi}_{\ell}(D)u_{\ell}=u_{\ell}$ and $h\equiv 1$ on
$\supp(g)$,
we have  $L_0+L_1=L$.
By Proposition ~\ref{cpct0}, it is enough to show the following three bounds:
First, there is $C(h)$, which only depends on $\max_{0\le |\alpha|\le [ r]+1}
\sup |\partial^\alpha h|$, so that
\begin{equation}\label{postcomp}
\|M  u\|_{C^ {p}_*}\le C(h)  \|u\|_{C^{p}_*}\, ,
\end{equation}
second, there is $C$, which does not depend on $\TT$ and $\gg$, so that  
for each $u\in C^p(K')$
\[
\|L'_0 u\|_{C^ {p}_*} \le
C \|\TT\|_+^ p \|\gg\|_{L^\infty} \|u\|_{C^{p}_*}\, 
,
\]
and finally,  for each $0<p'<p$ there is $C(\TT, \gg)$ so that
for each $u\in C^p(K')$
\begin{equation}\label{remember}
\|L'_1 u\|_{C^ {p}_*}< C(\TT,\gg)\|u\|_{C^{p'}_*}\, .
\end{equation}
(Note that if a Banach space $\cB'_1$ is compactly included in
a Banach space $\cB'_0$, then any bounded linear operator from
$\cB'_0$ to $\cB_1$ is compact when restricted to $\cB'_1$,
using that the composition of a compact operator followed by
a bounded operator is compact.)

Notice that there is $C$ (independent of $\TT$ and $\gg$) so that
\begin{equation}\label{cnormh}
\sum_{\ell: \ell\hookrightarrow n} 2^{pn-p\ell}
\le 2^{4p} \|\TT\|_+^p \sum_{j=0}^\infty  2^{-j}
\le C  \|\TT\|_+^p ,\, \forall n \, .
\end{equation}
Also notice that 
\begin{equation}\label{prop0}
\psi_m(D)\circ \psi_n(D)=(\psi_m\cdot \psi_n)(D)=0 \mbox{ when }
|m-n|\ge 5\, .
\end{equation}
The bound for $L'_0$ is then easy: 
\begin{align*}
\|L'_0 u\|_{C^ p_*} &=
\sup_m 2^{pm}
 \|\psi_m(D) \bigl (\sum_n( L_0 u)_{(n)}) \|_{L^\infty(\real^d)}\\
&\le \sup_m 2^{pm} \sum_{|n-m|<5}\sum_{\ell: \ell \hookrightarrow n} 
\|\psi_n(D) (L u_\ell) \|_{L^\infty}\\
&\le C \|\gg\|_{L^\infty} \sup_m 2^{pm} \sum_{|n-m|<5}
\sum_{\ell: \ell \hookrightarrow n} 
\|  u_\ell \|_{L^\infty}\\
&\le C \|\gg\|_{L^\infty}\sup_m  \sum_{|n-m|<5}
\bigl (
\sum_{\ell: \ell \hookrightarrow n} 
2^{pn-p\ell} \bigr )\|  u \|_{C^p_*}\\
&\le C \|\gg\|_{L^\infty} \|\TT\|^p_+ \|  u \|_{C^p_*}\, .
\end{align*}
We used 
$\psi_n(D) f =\widehat \psi_n * f$, which implies
$\|\psi_n(D) f \|_{L^\infty} \le C \|f\|_{L^\infty}$ for all
$n$, by  Young's inequality  for $L^\infty$.

\smallskip 
We will have to work a little harder for $L'_1$.
Assume first that  $p\le r-1$. 
Then it is enough to prove that for each $f\in C^ \infty(\real^d)$ with
rapid decay, and all
$n$
\begin{equation}\label{eqn:s0}
 \|\psi_{n}(D) (L(\tilde{\psi}_{\ell}(D) f)) \|_{L^\infty}
\le C(\TT,\gg)2^{-(r-1)\max\{n,\ell\}}  \|f\|_{L^\infty}\,
\mbox{if $\ell\not\hookrightarrow n$.}
\end{equation}
Indeed, using (\ref{prop0}) as in the estimate for $L'_0$, 
the above bound implies that  
\begin{align}
\|L'_1 u\|_{C^{p}_*}&
\le C(\TT,\gg) \cdot 
\sup_{n} \left(\sum_{\ell: \ell\not\hookrightarrow n}
 2^{p n-p' \ell-(r-1)\max\{n,\ell\}}\right)\|u\|_{C^{p'}_*}\, ,
\end{align}
and the conditions $p\le r-1$ and $p'>0$ ensure that the supremum over $n$
of the sum over $\ell$ such that $\ell\not\hookrightarrow n$
above is finite (recall the footnote~ 4).

To show (\ref{eqn:s0}),
we note that
\[
(\check\psi_{n}(D) L\,\tilde{\psi}_{\ell}(D) f )(x)
=(2\pi)^{-2d}\int_{\real^ d} V_{n}^{\ell}(x,y) \cdot f\circ \TT(y) |\det D\TT(y)|
dy \, ,
\]
where we have extended $\TT$ to a bilipschitz $C^r$ diffeomorphism
of $\real^d$ and 
\begin{equation}\label{Vkernel0}
V_{n}^{\ell}(x,y)=\int_{\real^d \times \real^d\times\real^d} 
e^{i(x-w)\xi+i(\TT(w)-\TT(y))\eta} 
\gg(w)\psi_{n}(\xi)
\tilde{\psi}_{\ell}(\eta)dw d\xi d\eta\, .
\end{equation}
Since $\| f\circ \TT \cdot |\det D\TT|\|_{L^\infty}\le C(\TT)\|f\|_{L^\infty}$, 
the inequality (\ref{eqn:s0}) 
follows if we show that  there exists
$C(\TT,\gg)$ such that for all $\ell\not\hookrightarrow n$
the operator norm of the integral operator 
\[
H^{\ell}_{n}:f\mapsto \int_{\real^d} V_{n}^{\ell}(x,y)f(y) dy
\]
acting on $L^\infty(\real^d)$ is 
bounded by $C(\TT,\gg)\cdot 2^{-(r-1)\max\{n,\ell\}}$. 

Define the integrable function $b:\real^d\to 
\real_+$ by
\begin{equation}\label{convol0}
b(x)=
1\quad\mbox{ if $\|x\|\le 1$}, \qquad
b(x)=\|x\|^{-d-1}\quad\mbox{ if $\|x\|> 1$.}
\end{equation}
The required estimate on $H^{\ell}_{n}$ follows if we show
\begin{equation}\label{eqn:Kernelest0}
|V_{n}^{\ell}(x,y)|\le C(\TT,\gg) 2^{-(r-1)\max\{n,\ell\}}\cdot  
2^{d\min\{n,\ell\}}b(2^{\min\{n,\ell\}}(x-y))\, ,
\end{equation}
for some  $C(\TT,\gg)>0$ and all  $\ell\not\hookrightarrow n$.
Indeed, as the right hand side 
of (\ref{eqn:Kernelest0}) is written as a function of $x-y$, say $B(x-y)$, we have, 
by Young's inequality  in $L^ \infty(\real^d)$,   
\begin{align*}
\|H^{\ell}_{n}f\|_{L^\infty}&\le \|B * f\|_{L^\infty}\le  \|B\|_{L^1} \|f\|_{L^\infty}\\
&\le C(\TT,\gg) 2^{-(r-1)\max\{n,\ell\}}\cdot \|b\|_{L^1}\cdot \|f\|_{L^\infty}\, .
\end{align*}
(Note that, by Young's inequality for $L^t(\real^d)$
with $1<t<\infty$,
the operator $H^{\ell}_{n}$
acting on each 
$L^ t(\real^ d)$  is also bounded by $C(\TT,\gg)\cdot 2^{-(r-1)\max\{n,\ell\}}$.
This is useful to control the essential spectral radius on anisotropic Sobolev
spaces, see \cite{BT}.)

We now prove  (\ref{eqn:Kernelest0}).
If $r\ge 2$ (otherwise we do nothing at this stage),
integrating (\ref{Vkernel0}) by parts $[r]-1$ times on $w$
(recall (\ref{parts})), we obtain 
\begin{equation}\label{firststep0}
V_{n}^{\ell}(x,y)=
\int e^{i(x-w)\xi+i(\TT(w)-\TT(y))\eta} 
F(\xi,\eta,w)\psi_{n}(\xi)
\tilde{\psi}_{\ell}(\eta)dwd\xi d\eta \, ,
\end{equation}
where $F(\xi,\eta,w)$ is a $C^{r-[r]}$ function in $w$
which is $C^\infty$ in the variables  $\xi$ and $\eta$.
The following exercise is an important (but straightforward) step
in the proof:

\begin{exercise}\label{forlater} 
Using (\ref{lowerbd0}), check 
that if $\psi_n(\xi)\cdot \tilde \psi_\ell(\eta)\ne 0$ then
\begin{equation}\label{firststep}
\|
F(\xi,\eta,\cdot)\|_{C^{r-[r]}}\le C(\TT,\gg)
2^{-([r]-1)\max\{n,\ell\}}\, .
\end{equation}
\end{exercise}

The estimate (\ref{firststep})
looks promising, but applying it naively  is not enough: since we are integrating over $\xi$
in the support of $\psi_n$ and over $\eta$ in the support of $\tilde \psi_\ell$,
we would get an additional factor $2^{dn+d\ell}$. In order to get rid of this factor,
we shall use another exercise:

\begin{exercise}
Using (\ref{lowerbd0}), show that if
$\psi_n(\xi)\cdot \tilde \psi_\ell(\eta)\ne 0$, then
for all multi-indices $\alpha$ and $\beta$
\begin{equation}\label{simplediff}
\|\partial_\xi^\alpha\partial_\eta^\beta
F(\xi,\eta,\cdot)\|_{C^{r-[r]}}\le C_{\alpha,\beta}(\TT,\gg)
2^{-n|\alpha|-\ell|\beta|-([r]-1)\max\{n,\ell\}}\, .
\end{equation} 
\end{exercise}

Assume first that $r$ is an integer (then, $r=[r]\ge 2$).
Put 
$$
G_{n,\ell}(\xi,\eta,w)=F(\xi,\eta,w)\psi_{n}(\xi)\tilde{\psi}_{\ell}(\eta)
$$ and
consider the scaling
$
\widetilde G_{n,\ell}(\xi,\eta,w)=G_{n,\ell}(2^n\xi, 2^\ell\eta,w)
$.

The estimate (\ref{simplediff}) implies that for all $\alpha$ and $\beta$
\begin{equation}\label{1more}
\|\partial_\xi^\alpha\partial_\eta^\beta
\widetilde G_{n,\ell}(\xi,\eta,\cdot)\|_{C^{r-[r]}}\le C_{\alpha,\beta}(\TT,\gg)
 2^{-([r]-1)\max\{n,\ell\}}\, ,
\forall \xi, \eta, n, \ell\, .
\end{equation}
Then, denoting by $\FF$ the inverse Fourier transform 
with respect to the variable $(\xi,\eta)$, and setting $W_{n}^{\ell}(u,v,w):=$
\begin{equation}\label{FFF}
(\FF \widetilde G_{n,\ell})(u,v,w)=(2\pi)^{-2d}\int_{\real^d}\int_{\real^d}
e^{iu\xi}e^{iv\eta} \widetilde G_{n,\ell}(\xi,\eta,w)\,  d\xi d\eta \, ,
\end{equation}
the bounds (\ref{1more})  imply   
that for any  nonnegative integers $k$ and $k''$ 
\begin{equation}
\label{rapiddecay}
\bigl \| \|u\|^k \|v\|^{k''}
W_{n}^{\ell}(u,v,\cdot) \bigr \|_{L^{\infty}}\le 
\widetilde C_{k,k'}(\TT,\gg) 2^{-([r]-1)\max\{n,\ell\}}\, ,
\forall u, v, n, \ell \, .
\end{equation}
(Just note that
the integrand in (\ref{FFF}) is supported in 
$\max\{\|\xi\|,\|\eta\|\}\le 2$, and integrate by parts with 
respect to $\xi$ and $\eta$ as many times as desired.)
Applying (\ref{rapiddecay}) to $k$, $k'$ in 
$\{0, d+1\}$, we get $C(\TT,\gg)$ so that for each $w\in K$, and
all $n$, $\ell$, $u$, $v$
\begin{equation}\label{bb}
|W_n^\ell(u,v,w)|\le C(\TT,\gg) 2^{-([r]-1)\max\{n,\ell\}}
b(u) b(v)  \, .
\end{equation}
(For $w\notin K$ we have $W_n^\ell(u,v,w)=0$ for all $u$, $v$, $n$, $\ell$.)
Therefore, since 
$$
(\FF G_{n,\ell})(u,v,w)=2^{dn+d\ell } W_{n}^{\ell}(2^{n}u,2^{\ell}v,w)\, ,
$$
we get by definition, 
\begin{align*}
&|V_{n}^{\ell}(x,y)|\le \int_K |(\FF G_{n,\ell})(x-w,\TT(w)-\TT(y),w) |\, dw\\
\nonumber &\, \le C \int_K 2^{dn +d\ell}
|W_{n}^{\ell}(2^{n}(x-w), 2^{\ell}(\TT(w)-\TT(y)),w)| \, dw \\
\nonumber &\, \le C(\TT,\gg)2^{-([r]-1)\max\{n,\ell\} +dn+d\ell}
\int_K    b(2^n(x-w))b(2^\ell(\TT(w)-\TT(y))  dw \, .
\end{align*}
Next,   using $u=2^n(x-w)$, note
$w_u=x-2^{-n}u$, and  write
\begin{align}
\label{above}
&\int_K 2^{dn +d\ell} b(2^n(x-w))b(2^\ell(\TT(w)-\TT(y))) \, dw\\
\nonumber &=\int_{\real^d}  2^{d\ell} b(u) b(2^\ell(\TT(w_u)-\TT(y)))\, du\, .
\end{align}
Since $\ell\le n+N(\TT)$ (see footnote 4),  we get
by using 
\begin{equation}\label{bu}
\int b(u)b(2^\ell(\TT(w_u)-\TT(y)))\, du\le \int b(u)\, du <\infty\, ,
\end{equation}
that
$|V_{n}^{\ell}(x,y)|\le C(\TT,\gg)2^{d\min\{n,\ell\}-([r]-1)\max\{n,\ell\}} $.

If $\|x-y\|> 2^{-\min\{n,\ell\}}$,  we can improve the estimate:
let $q_0\le \min\{\ell,n\}$ be the integer so that
$\|x-y\|\in [2^{-q_0}, 2^{-q_0+1})$.
Taking large constants $C(\TT)$, we may assume that for each $u\in \real^d$ either of the following conditions holds:
\begin{align*}
&\|u\|=2^{n}\|x-w_u\|\ge 2^{C(\TT)+n-q_0}\ge 2^{C(\TT)+\ell-q_0}\, \\
&2^{\ell}\|\TT(w_u)-\TT(y)\|>2^{C(\TT)+\ell}\|w_u-y\|>2^{C(\TT)+\ell-q_0}\, .
\end{align*}
Hence we obtain
\begin{align*}
\int &b(u)b(2^\ell(\TT(w_u)-\TT(y)))\, du \\
&\le
2^{C(\TT)+(n-q_0)(d+1)} \int  b(2^\ell(\TT(w_u)-\TT(y))) \, du\\
&\quad\qquad\qquad\qquad\quad +
2^{C(\TT)+(\ell-q_0)(d+1)} \int b(u)\, du\\
&\le  2^{C(\TT)-(d+1)(q_0-\ell)}.
\end{align*}
With this, we conclude
 $|V_{n}^{\ell}(x,y)|\le C(\TT,\gg) 2^{d\min\{n,\ell\}-([r]-1)\max\{n,\ell\}}2^{(d+1)q_0}$, proving (\ref{eqn:Kernelest0}) for integer $r$ .

If $r>1$ is not an integer, we 
start from (\ref{firststep0}) and rewrite $V_{n}^{\ell}(x,y)$
as  
\begin{equation}\label{firststep'}
\int e^{i\Lambda(x-w)(\xi/\Lambda)+i\Lambda(\TT(w)-\TT(y))(\eta/\Lambda)} 
F(\xi,\eta,w)\psi_{n}(\xi)
\tilde{\psi}_{\ell}(\eta)dwd\xi d\eta,
\end{equation}
for
$\Lambda=2^{\max\{\ell,n\}}$. Recalling (\ref{regparts}),
we apply to (\ref{firststep'})  one regularised integration by parts
for $\delta=r-[r]$ (noting that $\TT$ is $C^{1+\delta}$).
We get two terms $F_{1,\epsilon}(\xi,\eta,w)$
and $F_{2,\epsilon}(\xi,\eta,w)$.
Choosing   $\epsilon=\Lambda^{-1}$, we may apply the above procedure to each
of them.
The proof of (\ref{eqn:Kernelest0}) when $\gg$ is $C^ r$ and
$r-1< p\le r$ is done in Appendix \ref{appa0}.

It only remains to check (\ref{postcomp}) for the multiplication
operator $Mu=hu$.  Since $\TT$ is replaced by the identity map which satisfies
$\|\id\|_+\le 1$, and $\gg$ is replaced by
$h$, this can be done by a simplification
of the above arguments for $L'_0$ and
$L'_1$ (we can take $p'=p$),  decomposing
$M=M_0+M_1$ according to the relation $\hookrightarrow$ 
associated to $\id$. We leave details as an
exercise for the reader.
\end{proof}

\section{Bounding the essential spectral radius in the Anosov case}
\label{S1}

We now move to hyperbolic situations. We take $r>1$ and $X$ a compact
Riemann manifold, and assume that $T:X\to X$ is a $C^r$ Anosov diffeomorphism.
Recall that this means that $\Lambda=X$ is a hyperbolic set for $T$
in the sense of Section~\ref{state}.
For $g\in C^0(X)$, set
$$
R(T,g)= 
 \lim_{m \to \infty} (\sup  |g^{(m)}(x)|)^ {1/m}
$$
(the limit is well-defined and equal to the infimum,
by a standard subbadditivity argument).
We shall give the key steps
in the  proof of the following result
(which is weaker than Theorem~\ref{mainprop0}):

\begin{theorem}
[Essential spectral radius] \label{PF}
Let $T:X\to X$ be a $C^r$ Anosov diffeomorphism, and
let $g$ be a $C^ {r-1}$ function, 
 with $r>1$. For all real numbers $q<0<p$ with $p-q<r-1$, 
there is a
Banach space $C^{p,q}_*(T)$ of distributions on
$X$, containing all $C^s$ functions with $s> p$, and contained
in the dual of $C^s(X)$ for all $s<|q|$, on
which $\L_{T,g}$ extends boundedly and so that
\[
r_{ess}(\L_{T,g}|_{C^{p,q}_*(T)})\le
R(T,g) \max\{\lambda^p_s, \nu_u^q\}\, .
\]
In particular the pull-back operator $T^*\varphi=
\varphi \circ T$ satisfies for all $q<0<p$ with
$p-q<r-1$
$$
 r_{ess}(T^*|_{C^{p,q}_*(T)})\le
 \max\{\lambda^p_s, \nu_u^q\}<1\, .
$$
\end{theorem}

In the rest of this paper,
we explain how to adapt the tools in Section~\ref{toy} to prove the
above theorem.

\subsection{Local definition of the anisotropic norms.}
\label{SS1}

In this subsection we define the anisotropic norms in a compact
domain  of $\real^ d$. (The Banach space in Theorem~ \ref{PF} will be
constructed by patching together such local spaces in coordinate charts.)
Let  $\cone_+$ and $\cone_-$ be
closed cones in $\real^d$  with nonempty interiors, 
such that $\cone_+\cap\cone_-=\{0\}$. 
Let then
$\varphi_+, \varphi_-:\sphere\to [0,1]$ be
$C^{\infty}$ functions on the unit sphere $\sphere$ in $\real^d$ satisfying
\begin{equation}\label{vp}
\varphi_+(\xi)=
\begin{cases}
1, &\mbox{if $\xi\in \sphere\cap \cone_{+}$,}\\
0, &\mbox{if $\xi\in \sphere\cap \cone_{-}$,}
\end{cases} \qquad 
\varphi_-(\xi)=1-\varphi_+(\xi) \, .
\end{equation}
(What the reader can have in mind is that $\cone_+$ is a cone containing
a stable bundle and $\cone_-$ a cone containing an unstable bundle.

Except in Exercise~\ref{main3}, Theorem~\ref{th:main3}, and Exercise~\ref{postmult2}
below, we shall work in this
subsection with a fixed pair of cones $\cone_\pm$ and fixed functions
$\varphi_\pm$, they will not appear in the notation for the sake of
simplicity.  Recall $\psi_n$ and $\chi$ from Subsection~\ref{deff}.
For $n\in \mathbb{Z_+}$ and $\sigma\in\{+,-\}$,  we define
\[
\psi_{n, \sigma}(\xi)=
\begin{cases}
\psi_n(\xi)\varphi_{\sigma}(\xi/\|\xi\|),\quad&\mbox{ if $n\ge 1$,}\\
\chi(\|\xi\|)/2,\quad&\mbox{ if $n=0$.}
\end{cases}
\]

\begin{exercise} Prove that the $\psi_{ n, \sigma}$
enjoy similar properties as those of the $\psi_n$,
in particular the $L^1$-norm of the rapidly decaying function
$\widehat \psi_{n,\sigma}$ is bounded uniformly in $n$.
\end{exercise}

Fix $K\subset \real^d$ compact and with nonempty interior.
For  $u\in C^\infty(K)$,
define for each $n\in \mathbb Z_+$, $\sigma\in \{+,-\}$, and $x\in \real^d$:
\[
u_{n,\sigma}(x)=(\psi_{n,\sigma}(D) u)(x)
=(\widehat \psi_{n,\sigma} * u)(x)  \, .
\]
Since $1=\sum_{n= 0}^{\infty}\sum_{\sigma=\pm}\psi_{ n, \sigma}(\xi)$, 
we have $
u=\sum_{n\ge 0}\sum_{\sigma=\pm}u_{n,\sigma}$.

\begin{definition}[Anisotropic h\"older spaces $C^{p,q}_*(K)$]
Let $\cone_\pm$ and $\varphi_\pm$ be fixed, as above.
Let $p$ and $q$ be arbitrary real numbers.  
Define the anisotropic h\"older norm $\|u \|\hol{p,q}$
for   $u\in C^\infty(K)$,   by
\begin{equation}
\|u\|\hol{p,q}=\max\left\{\; 
\sup_{n\ge 0}\; 2^{pn} \| u_{n,+} \|_{L^\infty}, 
\; \; \sup_{n\ge 0} \; 2^{qn} \|u_{n,-} \|_{L^\infty}\;\right\}\, .
\end{equation}
Let $C^{p,q}_*(K)$  
be the completion of $C^{\infty}(K)$ for the norm 
$\|\cdot \|\hol{ p,q}$.  
\end{definition}

\begin{remark}
\label{justif}
In our application, $p>0$, and $q<0$. Recalling Section~\ref{toy}
for contracting branches and $p>0$,
it is then natural that $\cone_+$ and $\varphi_+$ be associated to a
contracting (i.e., stable) cone for the dynamics 
and $\cone_-$ and $\varphi_-$ be associated to an
expanding (i.e., unstable) cone.
Elements of  $C^{p,q}_*(K)$   are distributions which are
at  least $p$-smooth in the directions in $\cone_+$ and at most
$q$-``rough" in the directions of $\cone_-$.
\end{remark}

We shall not give a proof of the following result, referring instead
to \cite{BT}. (The proof is 
based on the Ascoli-Arzel\`a theorem.)

\begin{proposition}  [Compact embeddings]\label{prop:cpt}
If $p'< p$ and $q'< q$, the inclusion
$C_*^{p,q}(K)\subset C_*^{p',q'}(K)$ is compact. 
\end{proposition}


\subsection{Compact approximation for local  hyperbolic maps}
\label{S2}

Let $r >1$.  Let $K, K'\subset \real^d$ be compact subsets with non-empty interiors, and take 
a compact  neighborhood $W$ 
of $K$.  Let $\TT:W\to K'$ be a $C^{r}$ diffeomorphism onto
its image.
Let $\gg:\real^d\to \C$ be a $C^{r-1}$ function such that $\supp(\gg)$
is contained in the interior of $K$. 
In this section we study the transfer operator
\[
L: C^{r-1}(K')\to C^{r-1}(K), 
\qquad L u(x)= \gg(x)\cdot u\circ \TT(x) \, .
\]

For a pair of cones $\cone_\pm$ as in Subsection \ref{SS1},
we make   the following {\it cone-hyperbolicity} assumption on $\TT$:
\begin{equation}\label{conehyp}
D\TT_x^{tr}(\real^d \setminus \mbox{interior}\,(\cone_{+}))
\subset \mbox{interior}\,(\cone_{-})\cup\{0\}\quad \mbox{for all $x\in W$,}
\end{equation}
where $D\TT_x^{tr}$ denotes the transpose of the derivative of $\TT$ at $x$. 
(The above condition is sufficient in the
neighbourhood of a hyperbolic fixed point.
More generally, it will be useful to allow more flexibility and to work
with two pairs of cones. See Exercise~\ref{main3} below.)

Put
\begin{align*}
&\|\TT\|_+=\sup_{x}
\sup_{0\neq D\TT_x^{tr}(\xi)\notin \cone_{-}} \frac{\|D\TT_x^{tr}(\xi)\|}
{\|\xi\|}
\quad\hbox{(the ``weakest contraction")}\, , \\
&\|\TT\|_-=\inf_{x}
\inf_{0\neq \xi \notin \cone_{+}} \frac{\|D\TT_x^{tr}(\xi)\|}{\|\xi\|}\quad \hbox{(the ``weakest expansion ")} \, .
\end{align*}

The following result is the key to the proof of Theorem~\ref{PF}
(see also Exercise~\ref{main3} and Theorem~\ref{th:main3} below):

\begin{theorem}
[Estimates for local cone-hyperbolic maps]\label{th:main2}
For any $q'<q<0<p'<p$ such that $p-q'<r-1$,  
there exists a constant $C$  so that for each $C^r$ diffeomorphism 
$\TT$ and each  $C^{r-1}$ function $\gg$  as above 
(assuming in particular (\ref{conehyp})),
there is
a linear operator $L_1'$
such that for any  $u\in C_*^{p,q}(K')$
\begin{equation}\label{ess}
\|L u -L'_1 u\|_{C_*^{p,q}}\le C\|\gg\|_{L^\infty}\cdot 
 \max\{   \|\TT\|_+^{p},\|\TT\|_-^{q} \}\|u\|_{C_*^{p,q}}\, ,
\end{equation}
and, in addition, there is
$C(\TT, \gg)$ such that for any $u\in C_*^{p',q'}(K')$
$$
\|L'_1u\|_{ C^{p,q}_*}\le C(\TT,\gg)
\|u\|_{C^{p',q'}_*}\, .
$$ 

\end{theorem}

It is essential that the constant $C$ in (\ref{ess}) does not depend on
$\TT$ and $\gg$.

\begin{proof}[Proof of Theorem \ref{th:main2}]
We need more notation. 
By (\ref{conehyp}) there exist a  closed cone
$\widetilde \cone_+$ contained in the interior of
$\cone_+$ such that for all $x\in W$
\begin{equation}\label{conehyp2}
D\TT^{tr}_x(\real^d\setminus \mbox{interior}(\widetilde \cone_{+})) 
\subset  \, \mbox{interior}(\cone_{-})\cup \{0\} \, .
\end{equation}

Fix also a closed cone $\widetilde \cone_-$ contained in the interior
of $\cone_-$ and let  $\tilde\varphi_\pm:\sphere\to [0,1]$ be 
$C^{\infty}$ functions satisfying
\[
\tilde\varphi_-(\xi)=
\begin{cases}
0, &\mbox{if $\xi\in \sphere\cap \widetilde\cone_{+}$,}\\
1, &\mbox{if $\xi\notin \sphere\cap \cone_{+}$,}
\end{cases}
\qquad 
\tilde\varphi_+(\xi)=
\begin{cases}
1, &\mbox{if $\xi\notin \sphere\cap\cone_{-}$,}\\
0, &\mbox{if $\xi\in \sphere\cap \widetilde  \cone_{-}$.}
\end{cases} 
\]
Recalling $\tilde \psi_\ell$ from the beginning
of the proof of Theorem \ref{th:main1}, define for $\sigma\in \{+,-\}$ 
\begin{align*}
&\tilde \psi_{\ell, \sigma}(\xi)=
\begin{cases}
\tilde{\psi}_\ell(\xi)\tilde \varphi_{\sigma}(\xi/\|\xi\|),& \mbox{ if 
$\ell\ge 1$,}\\
\chi(2^{-1}\|\xi\|),& \mbox{ if $\ell=0$.}
\end{cases}
\end{align*}
Note that $\tilde \psi_{\ell,\tau}(\xi)=1$ 
if $\xi \in \supp (\psi_{\ell,\tau})$.

Up to slightly changing the cones $\widetilde \cone_{\pm}$, we can guarantee that
\begin{align}\label{notmuchchange}
&\inf_{x\in K}
\inf_{0\neq \xi \notin \widetilde \cone_{+}} \frac{\|D\TT_x^{tr}(\xi)\|}{\|\xi\|}
\ge \|\TT\|_- /2 \, ,
\\
&\sup_{x\in K}
\sup_{0\neq D\TT^{tr}_x(\xi) \notin \widetilde \cone_{-}} \frac{\|D\TT_x^{tr}(\xi)\|}{\|\xi\|}
\le 2 \|\TT\|_+ \,  \, .
\end{align}

We write $(\ell, \tau) \hookrightarrow (n,\sigma)$ if either 
\begin{itemize}
\item $(\tau,\sigma)=(+,+)$ and $2^ n\le 2^{\ell+5}\|\TT\|_+$, or 
\item $(\tau,\sigma)=(-,-)$  and 
$2^ {\ell-5} \|\TT\|_- \le 2^ n$, or
\item $(\tau,\sigma)=(+,-)$ and 
$2^ {n}\ge  2^5\|\TT\|_{-}$ or $2^ {\ell} \ge 2^5\| \TT\|_{+}$.
\end{itemize}
We write $(\ell, \tau) \not\hookrightarrow (n,\sigma)$ otherwise. 

\begin{exercise} 
Let $\check \cone_\pm$ be two closed cones with disjoint
interiors, so that $\check \cone_+\cap \check \cone_-=\{0\}$,
and with
$$
\mbox{closure}\, (\real^d \setminus \cone_+) 
\subset \mbox{interior}\, (\check \cone_-) \cup \{0\} \, .
$$
Let  $\check\varphi_\pm:\sphere\to [0,1]$, and 
$\check \psi_{\ell, \sigma}:\real^d \to [0,1]$,  for $\sigma\in \{+,-\}$
and $\ell \in \integer_+$, be functions defined
just like $\varphi_\pm$ 
and $\psi_{\ell, \sigma}$, but replacing
the cones $\cone_\pm$ by $\check \cone_\pm$.
Using  (\ref{conehyp2}) and (\ref{notmuchchange}), check
that there exists  an integer $N(\TT)>0$ such that for all $x\in \supp(\gg)$ 
\begin{equation}\label{lowerbd}
d(\supp(\check \psi_{n,\sigma}), D\TT_{x}^{tr}(\supp(\tilde \psi_{\ell,\tau})))
\ge 2^{\max\{n,\ell\}-N(\TT)} 
\quad \mbox{if $(\ell, \tau) \not\hookrightarrow (n,\sigma)$.}
\end{equation}
Hint: For $(\tau,\sigma)=(-,+)$,
use  (\ref{conehyp2}).
See \cite{BT} for further details. 

Note that (\ref{lowerbd}) is {\it exactly} the same lower bound as
(\ref{lowerbd0}).
\end{exercise}

Define  $L'_1$ and $L'_0$
by $L'_j u=\sum_{n,\sigma} (L_j u)_{ (n,\sigma)}$ with 
\[
(L'_0 u)_{ (n,\sigma)}
=\sum_{(\ell,\tau): (\ell,\tau)\hookrightarrow (n,\sigma)}
\check\psi_{ n,\sigma}(D) (L \, u_{ \ell,\tau})\, ,
\]
and
\[
(L'_1 u)_{ (n,\sigma)}=
\sum_{(\ell,\tau): (\ell,\tau) \not\hookrightarrow (n,\sigma)}
\check\psi_{ n,\sigma}(D) (L\,\tilde{\psi}_{\ell,\tau}(D) u_{\ell,\tau}) \, .
\]
Since $\tilde{\psi}_{\ell,\tau}(D)u_{\ell,\tau}=u_{\ell,\tau}$, 
we have  $L'_0+L'_1=L$.
Note also that by definition
of the cones $\check \cone_\pm$, if $|n-m|>5$ or $\upsilon=+$ and $\sigma=-$ then
for $i=0$ and $i=1$:
\begin{equation}\label{prop}
 \psi_{m,\upsilon}(D)
(L_i u)_{(n,\sigma)} =0 \, .
\end{equation}
By Proposition ~\ref{prop:cpt}, it is enough to show that
there is $C$, which does not depend on $\TT$ and $\gg$, so that  
for each $u\in C^{p,q}_*(K)$
\[
\|L'_0 u\|_{C_*^{p,q}}<C \max\{\|\TT\|_+^ p,\|\TT\|_-^ q\} \|\gg\|_{L^\infty} 
\|u\|_{C^{p,q}_*}\, ,
\]
and that for each $0<p'<p$ and $q'<q$
so that $p-q'<r-1$, there is $C(\TT, \gg)$ so that
for each $u\in C^{p,q}_*(K)$
\[
\|L'_1 u\|_{C^{ p,q}_*}< C(\TT,\gg)\|u\|_{C^{ p',q'}_*}\, .
\]

The bound for $L'_0$ is easy, like in the proof of Theorem~\ref{th:main1}:
Notice  that 
there is $C$ so that, setting  $c(+)=p$, $c(-)=q$, 
\begin{equation}\label{cnorm}
\sum_{(\ell,\tau): (\ell,\tau)\hookrightarrow (n,\sigma)} 2^{c(\sigma)n-c(\tau)\ell}
\le C
\max\{\|\TT\|_+^ p,\|\TT\|_-^ q\} ,\, \forall (n,\sigma)\, ,
\end{equation}
and  recall that
$\sup_{(n,\sigma)}
\int |\widehat\psi_{n,\sigma}(x)| dx < \infty$.

\smallskip

Consider  next $L'_1$.  
It is enough to prove that if $(\ell,\tau)\not\hookrightarrow(n,\sigma)$
then
\begin{equation}\label{eqn:s}
\|
\check \psi_{ n,\sigma}(D) (L\,\tilde{\psi}_{\ell,\tau}(D) f)\|_{L^\infty}
\le C(\TT,\gg)2^{-(r-1)\max\{n,\ell\}}\|f\|_{L^\infty}\, .
\end{equation}
Indeed, setting $c'(+)=p'$, and $c'(-)=q'$, (\ref{eqn:s}) 
and (\ref{prop}) imply that  
\begin{align*}
&\|L'_1 u\|_{C^{p,q}_*}\\
&\,\le \sup_{(m,\upsilon)} 
\sum_{(n,\sigma)}
\sum_{(\ell,\tau):(\ell,\tau)\not\hookrightarrow(n,\sigma)}
2^{c(\upsilon) m}\| \psi_{ m,\upsilon}(D)
\check \psi_{ n,\sigma}(D) (L\,\tilde{\psi}_{\ell,\tau}(D)
 u_{\ell,\tau})\|_{L^\infty}
\\
&\,\le C(\TT,\gg) \cdot 
\sup_{(n,\sigma)} \left(\sum_{(\ell,\tau):(\ell,\tau)\not\hookrightarrow(n,\sigma)}
 2^{c(\sigma) n-c'(\tau) \ell-(r-1)\max\{n,\ell\}}\right)\|u\|_{C_*^{p',q'}}\, .
\end{align*}
Then, since $p\le r-1$,  $p-q'< r-1$, and thus 
$-q< r-1$, we see from
the definition of $\not\hookrightarrow$ that 
\begin{equation}\label{star}
\sup_{(n,\sigma)} \left(\sum_{(\ell,\tau):(\ell,\tau)\not\hookrightarrow(n,\sigma)}
 2^{c(\sigma) n-c'(\tau) \ell-(r-1)\max\{n,\ell\}}\right)< \infty \, .
\end{equation}
(Note that $p-q\le r-1$ is not enough to guarantee the above bound
because of the case $(\tau,\sigma)=(-,+)$.)

To show (\ref{eqn:s}),
extend $\TT$ to $\real^d$ as in the
proof of Theorem~\ref{th:main1}, and rewrite  
\[
(
\check \psi_{n,\sigma}(D) (L\,\tilde{\psi}_{\ell,\tau}(D) f)(x)
=(2\pi)^{-2d}\int V_{n,\sigma}^{\ell,\tau}(x,y) \cdot f\circ \TT(y) |\det D\TT(y)|
dy \, ,
\]
where
\begin{align}\label{Vkernel}
&V_{n,\sigma}^{\ell,\tau}(x,y)=
\int e^{i(x-w)\xi+i(\TT(w)-\TT(y))\eta} 
\gg(w) \check{\psi}_{n,\sigma}(\xi)
\tilde{\psi}_{\ell,\tau}(\eta)dw d\xi d\eta \, .
\end{align}
Recall $b$ from (\ref{convol0}).
If we show
\begin{equation}\label{eqn:Kernelest}
|V_{n,\sigma}^{\ell,\tau}(x,y)|\le C(\TT,\gg) 2^{-(r-1)\max\{n,\ell\}}\cdot  
2^{d\min\{n,\ell\}}b(2^{\min\{n,\ell\}}(x-y))\, ,
\end{equation}
for some  $C(\TT,\gg)>0$ and all  $(\ell,\tau)\not\hookrightarrow (n,\sigma)$
then (\ref{eqn:s}) follows from Young's inequality, as in the expanding case
from Section~\ref{comcon}.

Finally, the proof of (\ref{eqn:Kernelest})
is {\it exactly} the same as the proof of (\ref{eqn:Kernelest0}), up to 
using
the change of variable $v=2^{\ell}(\TT(w)-\TT(y))$ instead
of $u=2^n(x-w)$ in (\ref{above}) if $\ell > n$.
\end{proof}

\begin{exercise}\label{main3}
Consider now two pairs of cones
$\cone_\pm$ and $\cone_\pm'$, and construct, for each
$p$ and $q$, two norms $\|\cdot\|_{C_*^{p,q}}$ and $\|\cdot\|_{(C'_*)^{p,q}}$
(by choosing $\varphi_\pm$ and $\varphi'_\pm$ as above).
Introduce a more general condition for $\TT$:
\begin{equation}\label{conehyp3}
D\TT_x^{tr}(\real^d \setminus \mbox{interior}\,(\cone'_{+}))
\subset \mbox{interior}\,(\cone_{-})\cup\{0\}\quad \mbox{for all $x\in W$.}
\end{equation}

Put
\begin{align*}
&\|\TT\|_+=\sup_{x\in K}
\sup_{0\neq D\TT_x^{tr}(\xi)\notin \cone_{-}} \frac{\|D\TT_x^{tr}(\xi)\|}
{\|\xi\|}
\quad\hbox{(the ``weakest contraction")}\, , \\
&\|\TT\|_-=\inf_{x\in K}
\inf_{0\neq \xi \notin \cone'_{+}} \frac{\|D\TT_x^{tr}(\xi)\|}{\|\xi\|}
\quad \hbox{(the ``weakest expansion")} \, .
\end{align*}

Check that a small modification
of the proof of Theorem \ref{th:main2} gives:

\begin{theorem}\label{th:main3}
For any $q'<q<0<p'<p$ such that $p-q'<r-1$
there exist a constant $C$ so that for each $C^r$ diffeomorphism
$\TT$ and $C^{r-1}$ function $\gg$,  assuming  (\ref{conehyp3}), there
exists   a linear operator $L'_1$ such that for any $u\in (C'_*)^{p,q}(K')$
\begin{equation*}
\|L u -L'_1 u\|_{C_*^{p,q}}\le C\|\gg\|_{L^\infty}\cdot 
 \max\{   \|\TT\|_+^{p},\|\TT\|_-^{q} \}\|u\|_{(C'_*)^{p,q}}\, ,
\end{equation*}
and, in addition, there is $C(\TT,\gg)$
so that for any $u\in (C'_*)^{p',q'}(K')$
$$
\|L'_1 u\|_{C_*^{p,q}}\le C(\TT,\gg) \|u\|_{( C'_*)^{p',q'}}\, .
$$  
\end{theorem}
(See \cite{BT}.)
\end{exercise}

\begin{remark}\label{postmult2}
 Theorem~\ref{th:main3}  may be applied to
$\TT$ the identity map, i.e.,  the operator $M u=h\cdot u$ of multiplication by
a smooth function $h$, up to taking suitable pairs of cones in order
to guarantee cone-hyperbolicity.
\end{remark}

\subsection{Transfer operators for Anosov diffeomorphisms}
\label{S6}

We prove
Theorem~\ref{PF} by reducing  to the model of 
Subsections \ref{SS1} and \ref{S2}.

\begin{proof}[Proof of  Theorem~\ref{PF}]
We first define the  space  $C_*^{p,q}(T)$\ 
by using local charts
to patch the  anisotropic  spaces from
Subsection ~\ref{SS1}. Fix 
a finite system of $C^{\infty}$ local charts $\{(V_j, \kappa_j)\}_{j=1}^{J}$ 
that cover $X$, 
and a finite system of
pairs of closed cones 
\footnote{We regard $\cone_{j,\pm}$ as constant cone fields in the {\em cotangent} bundle 
$T^*\real^d$.}  
$\{(\cone_{j,+}, \cone_{j,-})\}_{j=1}^{J}$ in $\real^d$
with the 
properties that for all  $1\le j,k\le J$:

\renewcommand{\labelenumi}{(\alph{enumi})}
\begin{enumerate}
\item The closure of $\kappa_j(V_j)$ is a compact subset $K_j$ of $\real^d$.
\item  
$\cone_{j,+} \cap \cone_{j,-}=\{0\}$.
\item If $x\in V_j$, the cones 
$(D\kappa_j)^{*}(\cone_{j,+})$ and 
$(D\kappa_j)^{*}(\cone_{j,-})$ in the cotangent space 
contain the normal subspaces of $E^s(x)$ and $E^u(x)$, respectively. 
\item  If 
$T^{-1}(V_k)\cap V_j \neq \emptyset$, setting 
$U_{jk}=\kappa_j(T^{-1}(V_k)\cap V_j)$,
the map in charts 
$
\TT_{jk}:=\kappa_k\circ T\circ \kappa_j^{-1}:U_{jk} \to \real^d
$
enjoys the cone-hyperbolicity condition:
\begin{equation}\label{charthyp}
D\TT_{jk,x}^{tr}(\real^d\setminus \mbox{interior } (\cone_{k,+}))
\subset \mbox{interior} (\cone_{j,-} )\cup\{0\}, \quad
\forall x \in U_{jk} .
\end{equation}
\end{enumerate}
The fact that such systems of cones exist is standard for
Anosov maps, see e.g. \cite{KH}.

Choose $C^\infty$ functions $\varphi_j^+, \varphi_j^-: \sphere\to [0,1]$ 
for $1\le j\le J$ 
which satisfy (\ref{vp}) with $\cone_\pm=\cone_{j,\pm}$, 
as in Section~\ref{SS1}.  This defines for each $j$
a local space denoted $C^{p,q,j}_*$.
Choose finally a $C^{\infty}$ partition of the unity
$\{\phi_j\}$ subordinate to the covering $\{V_j\}_{j=1}^J$, that is, the support of 
each $\phi_j:X\to [0,1]$ is contained in the interior 
of $V_j$, and 
we have $\sum_{j=1}^{J} \phi_j\equiv 1$ on $X$.

\begin{definition}
We define the Banach spaces $C_*^{p,q}(T)$ to be 
the completion of $C^{\infty}(X)$ for the norm
\[
\|u\|_{C_*^{p,q}(T)}:=\max_{1\le j\le J} 
\|(\phi_j\cdot u)\circ \kappa_j^{-1}\|_{C_*^{ p,q,j}} \, .
\]  
\end{definition}

By definition,  $C_*^{p,q}(T)$  contains 
$C^s(X)$ for $s>p$. 
If $0\le p'<p$ and $q'<q$,  
Lemma ~\ref{prop:cpt} and a finite
diagonal argument over $\{1, \ldots, J\}$, 
imply that the inclusion
$C^{p,q}_*(T) \subset C^{p',q'}_*(T)$  is compact.

For $m\ge 1$ and $j$, $k$ so that
\[
V_{m, jk}:=T^{-m}(V_k)\cap V_j \ne \emptyset \, ,
\]
we may consider the map  in charts 
\[
\TT^m_{jk}=\kappa_k\circ T^m\circ \kappa_j^{-1}:
\kappa_j(V_{m,jk})\to \real^d \, .
\]
Note that (\ref{charthyp}) implies that
\begin{equation*}
(D\TT^m_{jk,x})^{tr}(\real^d\setminus \mbox{interior } (\cone_{k,+}))
\subset \mbox{interior} (\cone_{j,-} )\cup\{0\}\, , \,\,\,
\forall x \in \kappa_j(V_{m,jk}) \, .
\end{equation*}

Set
\[
 R_{m}=  \max_{j,k}\sup_{x\in \kappa_j(V_{m,jk})}  
 |g^{(m)}\circ \kappa_j^{-1}(x)|\cdot
 \max\{ \|\TT^m_{jk}\|_+^p, \|\TT^m_{jk}\|_-^{q} \} \, ,
\]
where 
\[
\|\TT_{jk}^m\|_+=\sup_{x \in \kappa_j(V_{m,jk})}
\sup\left\{ \frac{\|(D\TT_{jk}^m)_x^{tr}(\xi)\|}{\|\xi\|}\;;\;
 0\neq (D\TT_{jk}^m)_x^{tr}(\xi)\notin \cone_{j,-}\right\} \, ,
\]
and
\[
\|\TT_{jk}^m\|_-=\inf_{x \in \kappa_j(V_{m,jk})}
\inf\left\{ \frac{\|(D\TT_{jk}^m)_x^{tr}(\xi)\|}{\|\xi\|}\;;\; 
 0\neq \xi \notin \cone_{k,+}\right\} \, .
\]
A standard argument in uniformly hyperbolic
dynamics gives
$$\lim_{m \to \infty} (\|\TT_{jk}^m\|_+)^{1/m}\le\lambda_s\, ,
$$
and
$$\lim_{m \to \infty} (\|\TT_{jk}^m\|_-)^{1/m}\ge \nu_u\, .
$$
Therefore
\begin{equation}\label{eta} 
\limsup_{m \to \infty} (R_{m})^ {1/m}\le  R(T,g) 
\max\{\lambda_s^p, \nu_u^q\}\, .
\end{equation}

Since $p-q<r-1$, 
we can apply Theorem~\ref{th:main3} to 
$\TT^m_{jk}$ and $\gg_j=(\phi_j g^{(m)})\circ \kappa_j^{-1}$
to obtain $C$ so that, setting $L_{jk}^{(m)}u=\gg_j\cdot (u
(\phi_k \circ \kappa_k^{-1}))\circ \TT^m_{jk}$
for $u\in C^{p,q,k}_*(K_k)$,
\[
\|L_{jk}^{(m)}u
-(L_{jk}^{(m)})'_1 u\|_{C_*^{p,q,j}}
\le C R_{m}\cdot \|u\|_{C^{p,q,k}_*} \, ,
\quad \forall m \, .
\]
with $\|(L_{jk}^{(m)})'_1(u)\|_{C^{p,q,j}_*}\le C(\TT^m_{jk}, \gg_j)
\|u\|_{C^{p',q',k}_*} $.
Using Remark ~\ref{postmult2}
and postcomposition by the multiplication operator
$M_j u= h_j \cdot u$ where $h_j:\real^d\to \infty$
is $C^\infty$, supported in $K_j$ and  $h_j\equiv 1$
on the support of $\phi_j \circ \kappa_j^{-1}$,
similarly as in the last paragraph of the proof
of Theorem~\ref{th:main1} (details are left to the
reader), this implies the claimed upper bound 
for the essential
spectral radius of $\mathcal{L}_{T,g}$. 
\end{proof}

\begin{remark}
Though it is not  explicit in our notation,
choosing a different system of local charts, a different partition of unity,
or a different set of cones or functions $\varphi_\pm$, does not a 
priori give rise to equivalent norms.  
This is a little unpleasant, but does not
cause  problems. 
\end{remark}


\appendix

\section{Theorem~\ref{th:exp}  when both $T$ and
$g$ are $C^r$}
\label{appa0}

\begin{proof}
We only need to adapt
the estimate (\ref{remember}) on $L_1$ to the case when $\gg$ is $C^{r}$
and $r-1< p\le r$, for $r>1$, for some $0<p'<p$.
Recall $V_n^\ell$ from (\ref{Vkernel0})
and $b$ from (\ref{convol0}).
We shall show 
\begin{equation}\label{eqn:Kernelest00}
|V_{n}^{\ell}(x,y)|\le C(\TT,\gg) 2^{-r\max\{n,\ell\}}\cdot  
2^{(d+1)\min\{n,\ell\}}b(2^{ \min\{n,\ell\}}(x-y)),
\end{equation}
for some  $C(\TT,\gg)>0$ and all  $\ell\not\hookrightarrow n$. 

\begin{exercise}
Show that (\ref{eqn:Kernelest00}) combined with 
$$
\sup_{n} \left(\sum_{\ell: \ell\not\hookrightarrow n}
 2^{p n-p' \ell +\min\{n,\ell\}-r\max\{n,\ell\}}\right) < \infty\, ,
$$
gives the claim.
(Recall footnote 4 and take $p'>p$ very close to $p$.)
\end{exercise}

Define for each $y$ a $C^r$ function:
\[
A_y(w)= \TT(w)-\TT(y)-D\TT(y)(w-y) \,  .
\]
We may rewrite
(\ref{Vkernel0})  as
\[
V_{n}^{\ell}(x,y)=\int e^{i(x-w)\xi+i D\TT(y) (w-y)\eta} 
\bigl ( e^{iA_y(w)\eta} \gg(w) \bigr ) \psi_{n}(\xi)
\tilde{\psi}_{\ell}(\eta)dw d\xi d\eta \, .
\]

Integrating (\ref{Vkernel0}) by parts once on $w$, we obtain 
\begin{equation}\label{truc}
V_{n}^{\ell}(x,y)=
\int e^{i(x-w)\xi+i(\TT(w)-\TT(y))\eta} 
\check F(\xi,\eta,w)\psi_{n}(\xi)
\tilde{\psi}_{\ell}(\eta)dwd\xi d\eta \, ,
\end{equation}
where $\check F(\xi,\eta,w)$ is a $C^{r-1}$ function in $w$
which is $C^\infty$ in the variables  $\xi$ and $\eta$.
(We used properties of the derivative of an exponential
to ``reconstruct"  $e^{i(\TT(w)-\TT(y))\eta}$.)
Then, integrate (\ref{truc}) $[r]-1$ times by parts on $w$, giving
\begin{equation}\label{tildeF}
V_{n}^{\ell}(x,y)=
\int e^{i(x-w)\xi+i(\TT(w)-\TT(y))\eta} 
\widetilde F(\xi,\eta,w)\psi_{n}(\xi)
\tilde{\psi}_{\ell}(\eta)dwd\xi d\eta \, ,
\end{equation}
with
$\widetilde F(\xi,\eta,w)$  a $C^{r-[r]}$ function in $w$
which is  $C^\infty$ in the variables  $\xi$ and $\eta$.
By (\ref{lowerbd0}), if
$\psi_n(\xi)\cdot \tilde \psi_\ell(\eta)\ne 0$,
then we have for all $\alpha$ and $\beta$
\begin{equation}\label{simplediffr}
\|\partial_\xi^\alpha\partial_\eta^\beta
\widetilde F\|_{C^{r-[r]}}\le C_{\alpha,\beta}(\TT,\gg) 2^\ell
2^{-n|\alpha|-\ell|\beta|-[r]\max\{n,\ell\}} \, .
\end{equation}
(The price we have to pay for the first integration by parts is
the factor $2^\ell$. What we gained is $2^{-[r]\max\{n,\ell\}},$
with $[r]$ instead
of $[r]-1$.)
Then (\ref{simplediffr}) implies (\ref{eqn:Kernelest00}), just like 
in Section~\ref{comcon}
(recall that $\ell\le n$). 
\end{proof}



\bibliographystyle{amsplain}

\end{document}